Article

# Assignment of Freight Traffic in a Large-Scale Intermodal Network under Uncertainty[†]

**Majbah Uddin** [1,*], **Nathan N. Huynh** [2] and **Fahim Ahmed** [3]

[1] National Transportation Research Center, Oak Ridge National Laboratory, 1 Bethel Valley Road, Oak Ridge, TN, 37830, USA; E-Mail: uddinm@ornl.gov

[2] Department of Civil and Environmental Engineering, University of Nebraska-Lincoln, 2200 Vine St, 262D, Lincoln, NE 68583-0851, United States; E-Mail: nathan.huynh@unl.edu

[3] Associate System Performance Engineer, South Carolina Department of Transportation, 955 Park St, Columbia, SC 29202, United States; E-Mail: ahmedf@scdot.org

[†] This manuscript has been authored by UT-Battelle, LLC, under contract DE-AC05-00OR22725 with the US Department of Energy (DOE). The US government retains and the publisher, by accepting the article for publication, acknowledges that the US government retains a nonexclusive, paid-up, irrevocable, worldwide license to publish or reproduce the published form of this manuscript, or allow others to do so, for US government purposes. DOE will provide public access to these results of federally sponsored research in accordance with the DOE Public Access Plan (http://energy.gov/downloads/doe-public-access-plan).

* For correspondence.

## Abstract

This paper presents a methodology for freight traffic assignment in a large-scale road-rail intermodal network under uncertainty. Network uncertainties caused by natural disasters have dramatically increased in recent years. Several of these disasters (e.g., Hurricane Sandy, Mississippi River Flooding, and Hurricane Harvey) severely disrupted the U.S. freight transportation network, and consequently, the supply chain. To account for these network uncertainties, a stochastic freight traffic assignment model is formulated. An algorithmic framework, involving the sample average approximation and gradient projection algorithm, is proposed to solve this challenging problem. The developed methodology is tested on the U.S. intermodal network with freight flow data from the Freight Analysis Framework. The experiments consider three types of natural disasters that have different risks and impacts on transportation networks: earthquakes, hurricanes, and floods. It is found that for all disaster scenarios, freight ton-miles are higher compared to the base case without uncertainty. The increase in freight ton-miles is the highest under the flooding scenario; this is because there are more states in the flood-risk areas, and they are scattered throughout the U.S.

---

## 1. Introduction

Efficient management of freight movements is essential to support domestic e-commerce and international trade. Freight activities are directly related to a country's Gross Domestic Product and economic viability. In recent years, the U.S. transportation system has supported a growing volume of freight, and it is anticipated that this trend will continue in the coming years. For example, the U.S. transportation system moved a daily average of about 53.6 million tons of freight valued at more than $54 billion in 2021. Freight tonnage is projected to increase by about 43% between 2023 and 2050 [1]. To support the projected increase in freight volume, an efficient, reliable, and low-cost freight logistics system is necessary to keep the U.S. competitive in the global market.

Current freight forecasting methodologies assume that the transportation network is always functioning and is never disrupted (e.g., [2,3]). Hwang and Ouyang [2] provided a framework for freight train traffic assignment in a network where the network links (i.e., rail tracks) are always available. Uddin and Huynh [3] provided a methodology for road-rail freight traffic assignment in an intermodal network which considered that the network elements are never disrupted. The assumptions were made by the authors to simplify the scope of the studies and were appropriate for the problems addressed in those studies. Those studies did not consider the risks from weather-induced disruptions which have dramatically increased in recent years; several have occurred recently that severely affected the U.S. freight transportation network. The Mississippi River flooding impacted a major freight route, I-40 in Arkansas in 2011. Tropical storm Irene caused damage to over 5,000 miles of highways and 34 bridges in Vermont in 2011. Hurricane Sandy caused billions of dollars in damage and severely flooded streets and tunnels in the New York and New Jersey region in 2011 [4]. In 2018, flooding from Hurricane Florence caused the closure of more than 200 roads in South Carolina and more than 600 roads in North Carolina, including several stretches of I-95, which is a major freight route along the Eastern seaboard [5]. In 2022, the U.S. endured 18 separate weather-related disasters with losses exceeding $9 billion each, with a total cost of about $171 billion [6]. Given the growing occurrence of such disasters and their impact on the freight transportation network, there is a need to develop freight forecasting methods that address network uncertainties caused by natural disasters.

To this end, this paper proposes a stochastic model for the assignment of freight traffic, considering road, rail, and intermodal shipments, on a road-rail intermodal network that is subject to uncertainty. Intermodal transportation is defined as the use of at least two modes (e.g., road and rail) to move freight shipments in intermodal containers from an origin to a destination. Given the exact evaluation of the stochastic model is difficult, an algorithmic framework is proposed for solving the model. To account for uncertainties in a realistic manner, the U.S. natural disaster risk map [7] is used. The disaster types considered are earthquakes, hurricanes, and floods. For each disaster scenario, the model seeks an equilibrium assignment for a given set of freight traffic demands between origins and destinations and available modes (road-only, rail-only, and intermodal). A comparative analysis of different disaster scenarios is performed to assess their impacts on the resulting freight flows.

## 2. Literature Review

In this section, published literature is explored in terms of freight flows and routing models, and solution approaches that are particularly relevant to the methodology proposed in this paper.

### 2.1. Freight Assignment and Routing Models

Considering vulnerability of the transportation network, the models for freight flows and routing can be classified into two categories: models under normal conditions (i.e., without disruptions) and models under disruptions. In the following sub-sections, the scope and key characteristics of these mathematical models are discussed.

#### 2.1.1. Models Under Normal Conditions

Majority of researchers used mathematical programs that seek to minimize the cost of freight flows under normal conditions. These studies focused on freight flows in intermodal/multimodal networks [3,8,9], freight routing in rail [10] and intermodal network [11], and optimally locating intermodal terminals [12] or hub locations in a capacitated network [13]. Besides cost-minimization models, a few models seek to minimize the travel time of freight flow. Hwang and Ouyang [2] developed a model that minimizes the railroad travel time under user equilibrium taking into account shipper and carriers' route choices. Chang [14] formulated a multi-objective multimodal multicommodity flow problem with time windows and concave costs. A few researchers used simulation-assignment approach to capture the operational issues such as delays at terminals and yards, technological and policy changes [15]. Validation of this approach is achieved by applying it to a large-scale rail network [16].

Researchers also developed freight flow models to find investment priorities for new establishment and/or improvement of existing infrastructure in multimodal/intermodal networks [17,18]. In freight rail transportation network, researchers explored the level of service in terms of loaded and empty freight flow, yard operations, and capacity constraints [19]. Other researchers presented models to incorporate shipper-carrier roles and interactions. Friesz et al. [8] developed a model that considered the role of both shipper and carrier in freight assignment. Lastly, Agrawal and Ziliaskopoulous [20] used Variational Inequality to incorporate the market equilibrium from shippers' perspective.

#### 2.1.2. Models Under Disruptions

All the aforementioned studies assume that the freight transport network is always functioning and is never disrupted. Researchers developed mathematical formulation for freight assignment, routing, investment priorities, and network resilience incorporating disruption in their model to tackle network vulnerability. Routing-based models are formulated for uncapacitated and capacitated networks, survivability of flow network under multiple arc failures [21], rerouting of coal by rail under disruption [22], and routing intermodal freight under disruption scenarios [23], hazmat transportation with random yard disruptions [24], and routing with carbon dioxide emissions consideration [25].

Some researchers have studied the resilience of freight transportation networks. Resiliency of a network is quantified as the ability to recover from disruption by preventing, absorbing, and mitigating its effect [26]. The quantification of the resilience of intermodal freight transportation is implemented using an unsatisfied demand parameter [27] and pathway resistance function in a time-variant multi-scale network [28]. Decision models are presented to reroute freight flows when a delay on a link exceeds a tolerance limit [29], for pre-disruption preparation and post-disruption recovery activities within a budget [30], for reliable multi-commodity routing in an intermodal road-rail network [31], reliable rail intermodal transport network [32], and robust model for rail-truck intermodal transportation network with random disruptions at intermodal yards [33]. Frameworks are also developed for quantifying the resilience of intermodal freight networks under extreme events such as natural disasters including seismic events [34, 35].

### 2.2. Solution Approach to Models

The majority of mathematical programs with or without consideration of disruption used either exact or heuristic methods as their solution approach. Apart from these two methods, Evans algorithm [8], Gauss-Seidal Linear Approximation [9], Diagonalization [19], modified convex combination [2], and modified gradient projection algorithms [3] are used in solving assignment models under normal conditions. On the other hand, models under disruption consideration used improved Depth First Search [29], and Sample Average Approximation [23] for routing; Benders Decomposition, column generation, Monte Carlo simulation [27], and Integer L-shaped [30] algorithm for resilience models.

The above studies focused on achieving cost-effectiveness, time-saving, and operational efficiency under prevailing normal conditions. However, freight assignment and routing models should incorporate uncertainties due to natural and man-made disruptions. A few freight routing models address the uncertainties from disruptions but to date, no model has been developed to comprehensively assign freight flows in an intermodal freight network under equilibrium conditions with consideration for disruptions. This study seeks to fill this gap in the literature by (a) proposing a new model for freight flow assignment in rail-road intermodal networks under disruptions, (b) developing a framework to solve the proposed model, and (c) verifying and demonstrating the scalability of the model by applying it to an actual large-scale intermodal freight network.

## 3. Model Formulation

This study assumes that in the long run the activities carried out by shippers and carriers will be in equilibrium. That means the cost of any shipment cannot be lowered by changing mode, route, or both. Additionally, it is assumed that the cost on all used paths via different modes (road-only, rail-only, and intermodal) is equal for each origin-destination (O-D) pair and equal to or less than the cost on any unused path at equilibrium [36].

The model formulation assumes that a road-rail intermodal freight transportation network is represented by a directed graph $\mathcal{G} = (\mathcal{N}, \mathcal{A})$, where $\mathcal{N}$ is the set of nodes and $\mathcal{A}$ is the set of links joining them in the network. Set $\mathcal{N}$ consists of the set of freight zone centroid nodes $\mathcal{N}_c$, the set of road intersections $\mathcal{N}_t$, and the set of rail junctions $\mathcal{N}_l$, that is, $\mathcal{N} = \mathcal{N}_c \cup \mathcal{N}_t \cup \mathcal{N}_l$. Set $\mathcal{A}$ consists of the set of road segment $\mathcal{A}_t$, the set of rail tracks $\mathcal{A}_l$, and the set of terminal transfer links $\mathcal{A}_f$, that is, $\mathcal{A} = \mathcal{A}_t \cup \mathcal{A}_l \cup \mathcal{A}_f$. The road-rail intermodal terminals are modeled as network links. The flows are bi-directional on the terminal links. The end nodes of terminals have different nodes, that is, one from set $\mathcal{N}_t$ and the other from set $\mathcal{N}_l$. Origin and destination sets are represented by $\mathcal{O} \subseteq \mathcal{N}$ and $\mathcal{D} \subseteq \mathcal{N}$, respectively. Table 1 summarizes the notations used in the model and Figure 1 presents a hypothetical road-rail intermodal freight transportation network.

The capacity of each network link $a \in \mathcal{A}$ is disruption-scenario dependent, that is, capacities will be different depending on disruption-scenario sample $\xi \in \Xi$. A decision variable $x_{a\xi}$ is defined to represent the assigned freight flow on link $a \in \mathcal{A}$ under disruption-scenario sample $\xi \in \Xi$. Typically, rail tracks are shared by train in both directions. For that reason, the link delay on any rail track is dependent on the flow on it as well as the flow in the opposite rail track. In the following model, for train flow, $x_{a\xi}$ represents the flow from node $i \in \mathcal{N}_l$ to node $j \in \mathcal{N}_l$, and $x_{a'\xi}$ represents the flow from node $j \in \mathcal{N}_l$ to node $i \in \mathcal{N}_l$.

**Table 1.** Mathematical Notation.

| Notation | Description |
|---|---|
| $\mathcal{N}$ | set of nodes in network |
| $\mathcal{A}$ | set of links in network |
| $\mathcal{N}_c$ | set of freight zone centroid nodes in network |
| $\mathcal{N}_t$ | set of road intersections in network |
| $\mathcal{N}_l$ | set of rail junctions in network |
| $\mathcal{A}_t$ | set of road segments in network |
| $\mathcal{A}_l$ | set of rail tracks in network |
| $\mathcal{A}_f$ | set of terminal transfer links in network |
| $\mathcal{O}$ | set of origins in network, $\mathcal{O} \subseteq \mathcal{N}$ |
| $\mathcal{D}$ | set of destinations in network, $\mathcal{D} \subseteq \mathcal{N}$ |
| $T$ | set of available intermodal terminals for transfer of shipments |
| $o$ | origin zone index, $o \in \mathcal{O}$ |
| $d$ | destination zone index, $d \in \mathcal{D}$ |
| $K_t^{od}$ | set of paths with positive truck flow from $o$ to $d$ |
| $K_l^{od}$ | set of paths with positive train flow from $o$ to $d$ |
| $K_i^{od}$ | set of paths with positive intermodal flow from $o$ to $d$ |
| $q_t^{od}$ | freight truck demand from $o$ to $d$ |
| $q_l^{od}$ | freight train demand from $o$ to $d$ |
| $q_i^{od}$ | freight intermodal demand from $o$ to $d$ |
| $\Xi$ | set of disruption-scenario samples |
| $\xi$ | a disruption-scenario sample, $\xi \in \Xi$ |
| $f_{k\xi}^{od}$ | flow on path $k$ connecting $o$ and $d$ under disruption-scenario sample $\xi$ |
| $x_{a\xi}$ | flow on link $a \in \mathcal{A}$ under disruption-scenario sample $\xi$ |
| $C_{a\xi}$ | capacity of link $a \in \mathcal{A}$ under disruption-scenario sample $\xi$ |
| $t_{a\xi}(\omega)$ | travel time on link $a \in \mathcal{A}$ for flow of $\omega$ under disruption-scenario sample $\xi$ |

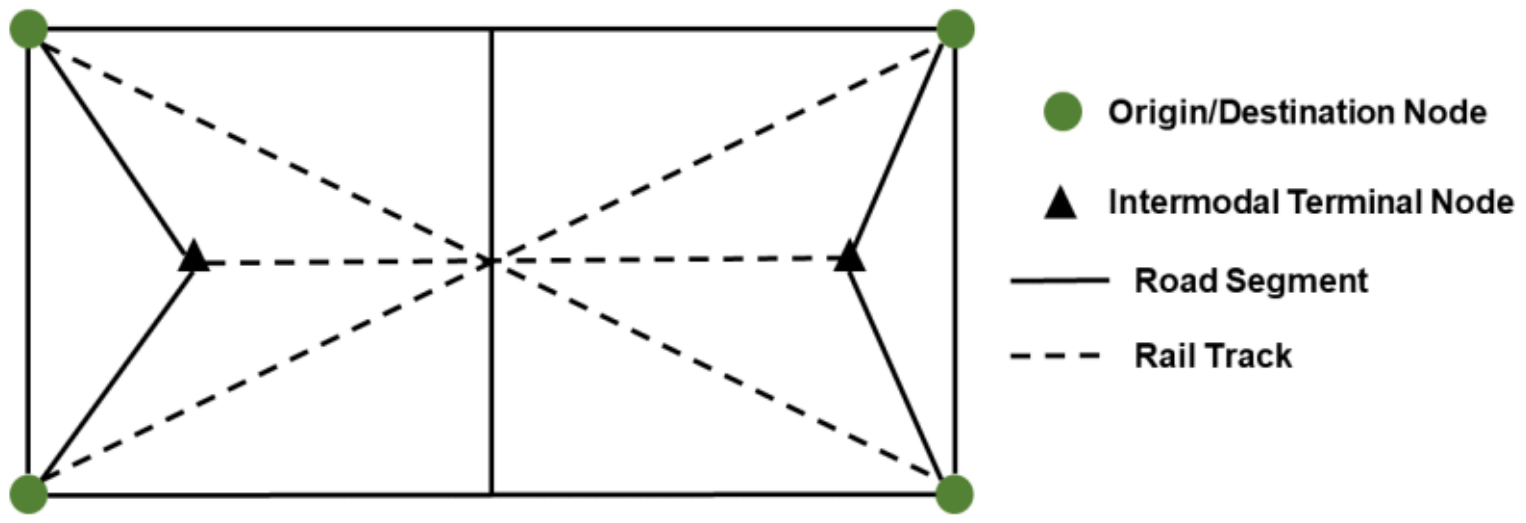

**Figure 1.** Hypothetical network.

For freight truck demand $q_t^{od}$ from origin $o \in \mathcal{O}$ to destination $d \in \mathcal{D}$ and a set of paths $K_t^{od}$ that connect $o$ to $d$ for each O-D pair, the path flow $f_{k\xi}^{od}$ satisfies the demand under disruption-scenario sample $\xi \in \Xi$ $\left(\sum_{k \in K_t^{od}} f_{k\xi}^{od} = q_t^{od}\right)$. Similarly, the path flows for freight train and intermodal on path sets $K_l^{od}$ and $K_i^{od}$ satisfy their respective demands ($q_l^{od}$ and $q_i^{od}$) from $o$ to $d$ under disruption-scenario sample $\xi$. Since the intermodal path set consists of paths formed by links from both road segments and rail tracks, the total freight flow on a road segment ($a \in \mathcal{A}_t$) is the sum of the road-only flows and the intermodal flows. Similarly, the total freight flow on a rail track ($a \in \mathcal{A}_l$) is the sum of the rail-only flows and intermodal flows. Using the parameters and decision variable described above, the following stochastic model finds the equilibrium freight flows in a road-rail intermodal network.

$$\text{Min} \quad \mathbb{E}_{\xi}\left[\sum_{a \in \mathcal{A}_t} \int_0^{x_{a\xi}} t_{a\xi}(\omega)\, d\omega + \sum_{a \in \mathcal{A}_l} \int_0^{x_{a\xi}+x_{a'\xi}} t_{a\xi}(\omega)\, d\omega\right] \tag{1}$$

Subject to

$$\sum_{k \in K_t^{od}} f_{k\xi}^{od} = q_t^{od}, \qquad \forall\, o \in \mathcal{O},\ d \in \mathcal{D},\ \xi \in \Xi \tag{2}$$

$$\sum_{k \in K_l^{od}} f_{k\xi}^{od} = q_l^{od}, \qquad \forall\, o \in \mathcal{O},\ d \in \mathcal{D},\ \xi \in \Xi \tag{3}$$

$$\sum_{k \in K_i^{od}} f_{k\xi}^{od} = q_i^{od}, \qquad \forall\, o \in \mathcal{O},\ d \in \mathcal{D},\ \xi \in \Xi \tag{4}$$

$$x_{a\xi} = \sum_{o \in \mathcal{O}} \sum_{d \in \mathcal{D}} \sum_{k \in K_t^{od}} f_{k\xi}^{od} \delta_{ka}^{od} + \sum_{o \in \mathcal{O}} \sum_{d \in \mathcal{D}} \sum_{k \in K_i^{od}} f_{k\xi}^{od} \delta_{ka}^{od}, \qquad \forall\, a \in \mathcal{A}_t,\ \xi \in \Xi \tag{5}$$

$$x_{a\xi} = \sum_{o\in\mathcal{O}} \sum_{d\in\mathcal{D}} \sum_{k\in K_l^{od}} f_{k\xi}^{od}\delta_{ka}^{od} + \sum_{o\in\mathcal{O}} \sum_{d\in\mathcal{D}} \sum_{k\in K_i^{od}} f_{k\xi}^{od}\delta_{ka}^{od}, \quad \forall\, a \in \mathcal{A}_l,\ \xi \in \Xi \tag{6}$$

$$f_{k\xi}^{od} \geq 0, \qquad \forall\, k \in K_t^{od},\ k \in K_l^{od},\ k \in K_i^{od},\ o \in \mathcal{O},\ d \in \mathcal{D},\ \xi \in \Xi \tag{7}$$

where

$$\delta_{ka}^{od} = \begin{cases} 1 & \text{if link } a \text{ is on path } k \text{ connecting } o \text{ and } d \\ 0 & \text{otherwise} \end{cases}$$

Constraints 2 through 4 ensure that all freight demands are assigned to the network. Constraints 5 and 6 compute the link flows on road and rail segments, respectively. Lastly, constraints 7 enforces all flow to be nonnegative. To estimate the objective function value in Equation 1, travel time on road and rail segments as a function of the flow are needed. For the road travel time, the Bureau of Public Roads (BPR) link performance function is used. The BPR function is considered given that it is most widely used. However, other performance functions could be used in conjunction with the model. The overall model output pattern would be the same for any specification for performance function. For rail travel time, the link performance function proposed by Uddin and Huynh [3] is used. The link performance functions have the following form:

$$t_{a\xi}(x_{a\xi}) = t_{0,t}\left(1 + 0.15\left(\frac{x_{a\xi}}{C_{a\xi}}\right)^4\right), \qquad \forall\, a \in \mathcal{A}_t,\ \xi \in \Xi \tag{8}$$

$$t_{a\xi}(x_{a\xi} + x_{a'\xi}) = t_{0,l}\left(1 + \left(\frac{x_{a\xi} + x_{a'\xi}}{C_{a\xi}}\right)^4\right), \qquad \forall\, a \in \mathcal{A}_l,\ \xi \in \Xi \tag{9}$$

$t_{0,t}$ and $t_{0,l}$ are the free-flow travel time for road and rail links, respectively. $C_{a\xi}$ is the capacity of the link $a$ under disruption-scenario sample $\xi$.

## 4. Algorithmic Strategy

The proposed model (1) – (7) is a stochastic program, which is difficult to solve because of the need to evaluate the expectation in the objective function. One approach is to approximate the expected value through sample averaging [23,37]. This approach is known as sample average approximation (SAA). In this study, the SAA algorithm proposed by Santoso et al. [37] is adopted. The SAA method is an approach for solving stochastic optimization problems by using Monte Carlo simulation. There are several advantages to using this method. The method can be applied to many stochastic optimization problems. It also has desirable convergence properties [38]. The objective function of the model (Equation 1) can be rewritten as follows, without loss of generality, where $y$ represents the decision variable.

$$\text{Min} \quad \mathbb{E}_{\xi}[Q(y,\xi)] \tag{10}$$

### 4.1. The SAA Algorithm

*Step 1*

Generate $M$ independent disruption-scenario samples each of size $N$, i.e., $\xi_m^1, \ldots, \xi_m^N$ for $m = 1, \ldots, M$. For each sample, solve the corresponding SAA problem.

$$\text{Min} \quad \frac{1}{N}\sum_{n=1}^{N} Q(y, \xi_m^n) \tag{11}$$

Let $z_N^m$ and $\hat{y}_N^m$, $m = 1, \ldots, M$, be the corresponding optimal objective value and an optimal solution, respectively.

*Step 2*

Compute the following two values.

$$\bar{z}_{N,M} := \frac{1}{M}\sum_{m=1}^{M} z_N^m \tag{12}$$

$$\sigma^2_{\bar{z}_{N,M}} := \frac{1}{(M-1)M} \sum_{m=1}^{M} \left( z_N^m - \bar{z}_{N,M} \right)^2 \tag{13}$$

The expected value of $z_N$ is less than or equal to the optimal value $z^*$ of the true problem [37]. Thus, $\bar{z}_{N,M}$ provides a lower statistical bound for the optimal value $z^*$ of the true problem, and $\sigma^2_{\bar{z}_{N,M}}$ is an estimate of the variance of this estimator.

*Step 3*

Choose a feasible solution $\tilde{y}$ from the above-computed solutions $\hat{y}_N^m$, and generate another $N'$ independent disruption-scenario sample, i.e., $\xi^1, \ldots, \xi^{N'}$. Then estimate true objective function value $\tilde{z}_{N'}(\tilde{y})$ and variance of this estimator as follows.

$$\tilde{z}_{N'}(\tilde{y}) := \frac{1}{N'} \sum_{n=1}^{N'} \mathcal{Q}(\tilde{y}, \xi^n) \tag{14}$$

$$\sigma^2_{N'}(\tilde{y}) := \frac{1}{(N'-1)N'} \sum_{n=1}^{N'} \left[ \mathcal{Q}(\tilde{y}, \xi^n) - \tilde{z}_{N'}(\tilde{y}) \right]^2 \tag{15}$$

Typically, $N'$ is much larger than the sample size $N$ used in solving the SAA problems. $\tilde{z}_{N'}(\tilde{y})$ is an unbiased estimator of $z(\tilde{y})$. Also, $\tilde{z}_{N'}(\tilde{y})$ is an estimate of the upper bound on $z^*$.

*Step 4*

Compute an estimate of the optimality gap of the solution $\tilde{y}$ using the lower bound estimate and the objective function value estimate from Steps 2 and 3, respectively, using the equations below:

$$\text{gap}_{N,M,N'}(\tilde{y}) := \tilde{z}_{N'}(\tilde{y}) - \bar{z}_{N,M} \tag{16}$$

The estimated variance of the above gap estimator is then given by

$$\sigma^2_{\text{gap}} = \sigma^2_{N'}(\tilde{y}) + \sigma^2_{\bar{z}_{N,M}} \tag{17}$$

4.2. Gradient Projection Algorithm

The SAA problem in Equation 11 is the standard traffic assignment problem, which cannot be solved analytically. This study adopts the path-based algorithm (gradient projection) proposed by Uddin and Huynh [3] to solve the traffic assignment problem. The gradient projection (GP) algorithm was first used by Jayakrishnan et al. [39] to solve the traffic assignment problem. Uddin and Huynh [3] further modified the GP algorithm to consider the situation where freight traffic demands could be transported via one of three modes (road-only, rail-only and intermodal). Their GP algorithm also considered intermodal terminals in the network. The adopted GP algorithm has the following iterative steps for a specific disruption-scenario sample $\xi$.

*Step 0 (Initialization)*

Set $t_{a\xi} = t_{a\xi}(0), \forall a \in \mathcal{A}$, and select terminals from the available terminals $T$ for all O-D pairs. Assign O-D demand $q_t^{od}$, $q_l^{od}$, and $q_i^{od}$ on the shortest path calculated based on $t_{a\xi}, \forall a \in \mathcal{A}_t$, $t_{a\xi}, \forall a \in \mathcal{A}_l$, and $t_{a\xi}, \forall a \in \mathcal{A}$, respectively, and initialize the path sets $K_t^{od}$, $K_l^{od}$, and $K_i^{od}$ with the corresponding shortest path for each O-D pair $(o, d)$. This initialization yields path flows and link flows. Set iteration count to $p = 1$.

*Step 1*

For each O-D pair $(o, d)$:

*Step 1.1 (Update)*

Set $t_{a\xi} = t_{a\xi}\left(x_{a\xi}(p)\right), \forall a \in \mathcal{A}$. Update the first derivative lengths, i.e., path travel times at current flow: $d_{kt}^{od}(p), \forall k \in K_t^{od}$, $d_{kl}^{od}(p), \forall k \in K_l^{od}$, and $d_{ki}^{od}(p), \forall k \in K_i^{od}$.

*Step 1.2 (Direction finding)*

Find the shortest path $\bar{k}_t^{od}(p)$ based on $t_{a\xi}(p), \forall a \in \mathcal{A}_t$. If different from all the paths in $K_t^{od}$, add it to $K_t^{od}$ and record $d^{od}_{\bar{k}_t^{od}(p)}$. If not, tag the shortest among the paths in $K_t^{od}$ as $\bar{k}_t^{od}(p)$. Repeat this procedure for $K_l^{od}$ and $K_i^{od}$ to find $d^{od}_{\bar{k}_l^{od}(p)}$ and $d^{od}_{\bar{k}_i^{od}(p)}$ based on $t_{a\xi}, \forall a \in \mathcal{A}_l$ and $t_{a\xi}, \forall a \in \mathcal{A}$, respectively.

*Step 1.3 (Move)*

Set the new path flows for $K_t^{od}$.

$$f_{k\xi}^{od}(p+1) = max\left\{0, f_{k\xi}^{od}(p) - \frac{\alpha(p)}{s_{k\xi}^{od}(p)}\left(d_{kt}^{od}(p) - d^{od}_{\bar{k}_t^{od}(p)}\right)\right\}, \forall k \in K_t^{od}, k \neq \bar{k}_t^{od} \tag{18}$$

where

$$s_{k\xi}^{od}(p) = \sum_a \frac{\partial t_{a\xi}^{od}(p)}{\partial x_{a\xi}^{od}(p)}, \forall k \in K_t^{od} \tag{19}$$

$a$ denotes the links that are on either $k$ or $\bar{k}_t^{od}$, but not on both. $\alpha(p)$ is the step size; the value of this parameter is set as 1 [39]. Now,

$$f^{od}_{\bar{k}^{od}\xi}(p+1) = q_t^{od} - \sum f_{k\xi}^{od}(p+1), \forall k \in K_t^{od}, k \neq \bar{k}_t^{od}(p) \tag{20}$$

Follow the above procedure to find new path flow for $K_l^{od}$ and $K_i^{od}$. From path flows find the link flows $x_{a\xi}(p+1)$.

*Step 2 (Convergence test)*

If the convergence criterion is met, stop. Else set $p = p + 1$ and go to Step 1.

Figure 2 shows a flow chart that illustrates how the SAA and GP algorithms are used to solve the traffic assignment problem. The model solution procedure starts with the input of O-D demands and intermodal network data. Then, a number of disruption-scenario samples are generated following the procedure described in Step 1 of the SAA algorithm. Then, for a specific scenario sample, GP algorithm solves the assignment problem and outputs the network link flows. This is repeated until all the scenario samples have been considered. After that, the procedure continues to Step 2 to 4 of the SAA algorithm.

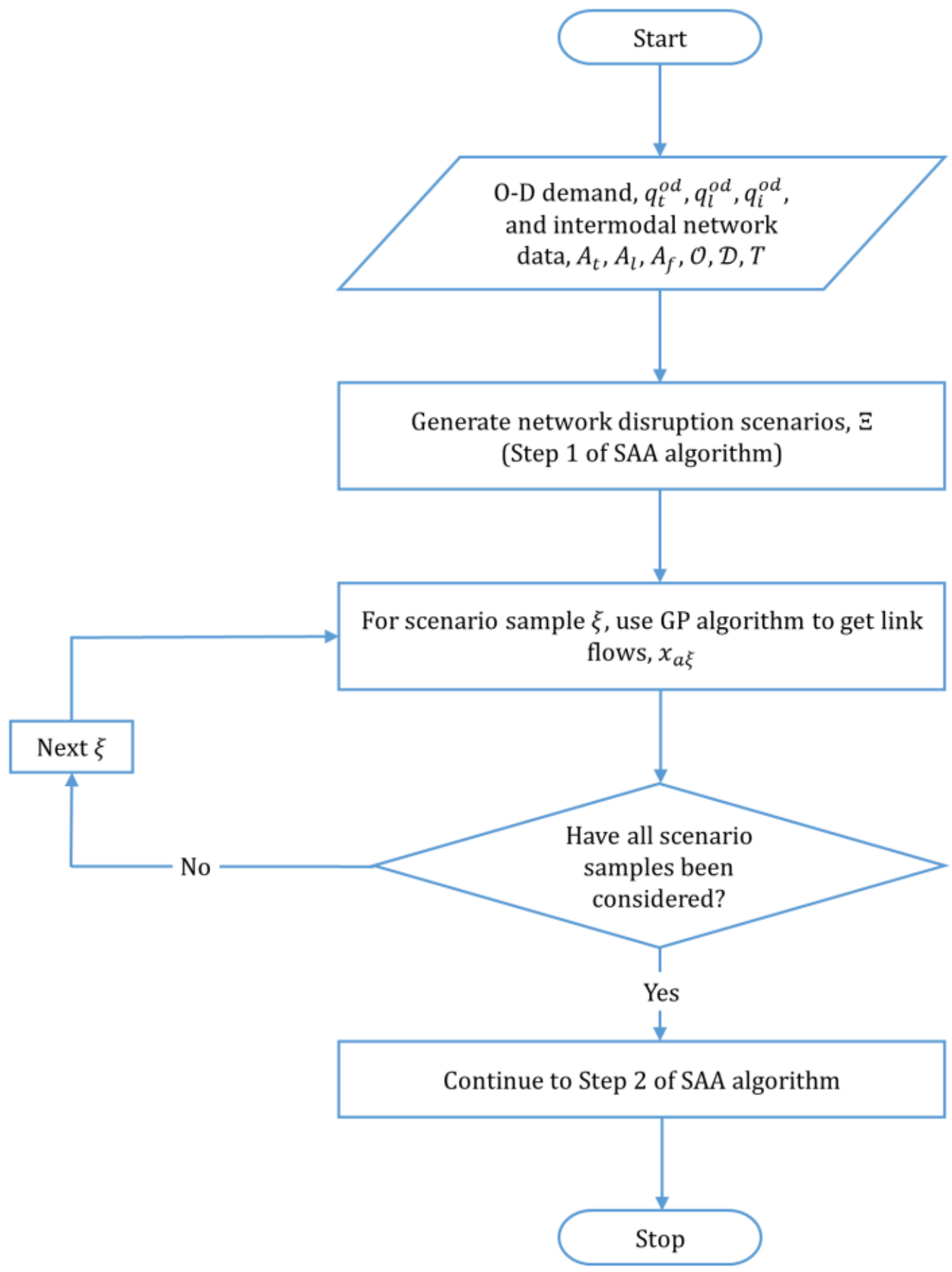

**Figure 2.** Algorithmic framework.

## 5. Numerical Experiments

The numerical experiments were run on a desktop computer with an Intel Core i7 3.40-GHz processor and 24 GB of RAM. To validate the proposed model and algorithmic framework, the road-rail transportation network in the contiguous U.S. and four disaster scenarios were considered.

### 5.1. Network and Disaster Data

The numerical experiments used the U.S. road-rail intermodal network shown in Figure 3 [3]. This network is a simplified version of the U.S. intermodal network created by Oak Ridge National Laboratory [40]. The simplified network consists of only Interstates, Class I railroads, and road-rail terminals. In Figure 3, the green circles represent Freight Analysis Zone (FAZ) centroids, the black triangles represent road-rail terminals, the black lines represent Interstates, and the grey dashed lines represent Class I railroads. In all, the network has a total of 1,532 links and 301 nodes. The nodes include 120 FAZ centroids, 97 major road intersections, and 84 major rail junctions.

The Freight Analysis Framework (FAF) is the most comprehensive public source of freight data in the U.S. [41]. Currently, FAF version 5 dataset is available. However, in this paper FAF version 3 (FAF3) was used due to two reasons: (i) the impacts of disaster scenarios were evaluated against the base case without uncertainty [3] which is based on FAF3 and (ii) the network used for experiments (Figure 3) were generated using inputs from FAF3. The FAF5 has 132 FAZ compared to 120 FAZ in FAF3. Note that the proposed model and algorithmic framework can assign freight flows using data from any version of FAF since FAZ is an input for O-D pairs.

It is assumed that there are 14,400 possible O-D demand pairs in the network since FAF3 has 120 FAZ. The freight demand from FAF is provided in terms of tonnage. Therefore, it needs to be converted to number of trucks or trains as input to the model. This study used the number of trucks and trains converted from freight demands from Uddin and Huynh [3]. The procedure to convert tonnage to truck is based on truck equivalency factor for different truck body types and empty truck factor. For intermodal freight transportation, half of the truck trips consist of empty movements (i.e., truck travels unloaded to its next destination) [42]. The procedure for train and intermodal demands are based on average loading capacity. The freight O-D trip tables for the truck, rail, and intermodal trips are in $120 \times 120 \times 3$ matrix form, corresponding to 120 origins, 120 destinations and 3 shipping modes. For a single day in the base year (2007), there are 618,190 truck shipments, 1,415 rail shipments, and 12,474 intermodal shipments. Note that FAF3 data are provided on an annual basis. For daily demand conversion, it was assumed that the road and rail infrastructure are operational for 365 days a year.

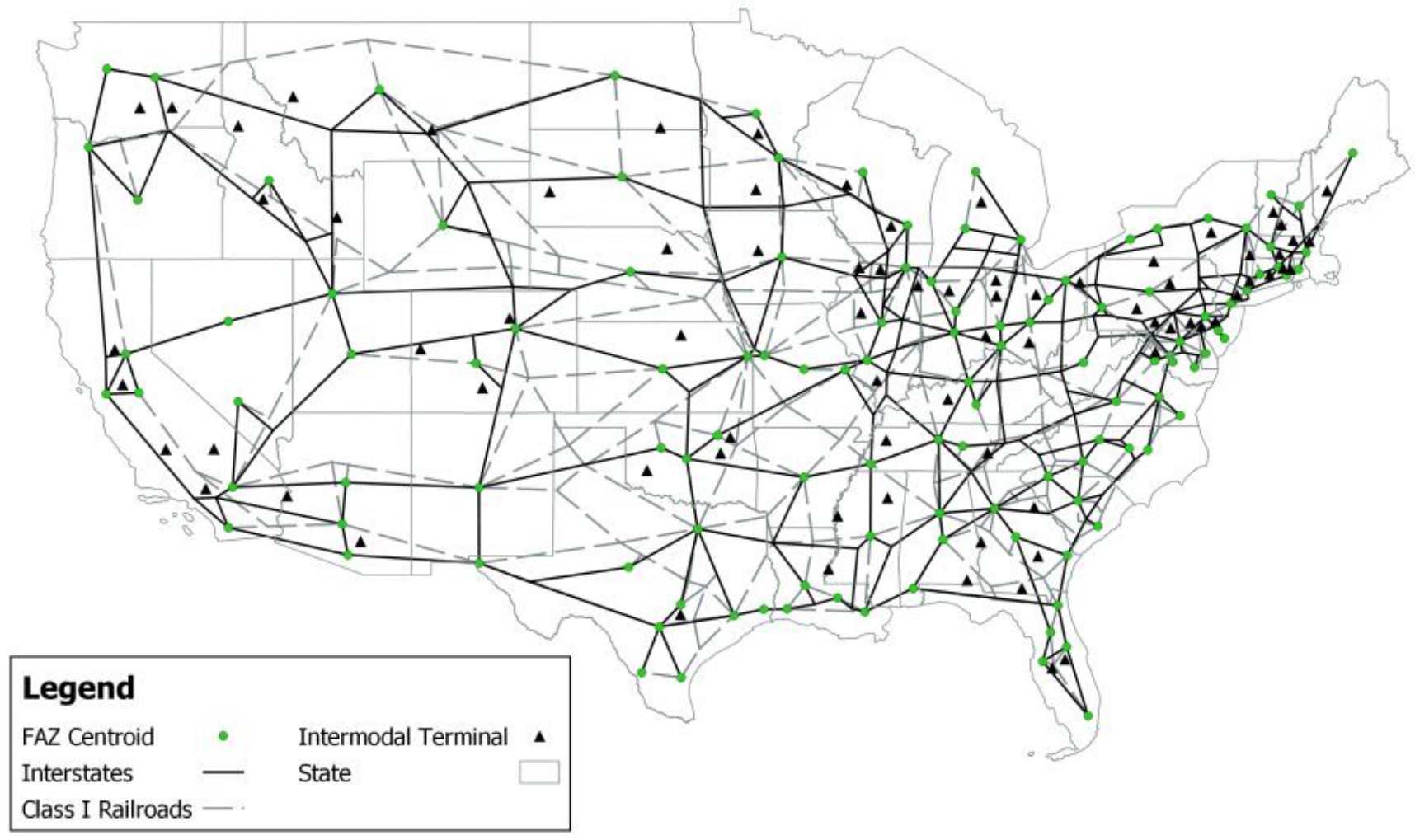

**Figure 3.** U.S. road-rail intermodal network [3].

To create the disaster scenarios, the U.S. natural disaster risk map shown in Figure 4 was used; this map was generated using data from the American Red Cross and the National Oceanic and Atmospheric Administration [7]. As shown, the high-risk earthquake is limited to California and a few other states, whereas the moderate-risk earthquake covers a significantly larger area. For this reason, two earthquake scenarios were considered: one with only high-risk areas, and another with both high and moderate risk areas. In all, four disaster scenarios were considered for the numerical experiments. The scenarios are earthquake (high risk), earthquake (high and moderate risk), hurricane, and flood.

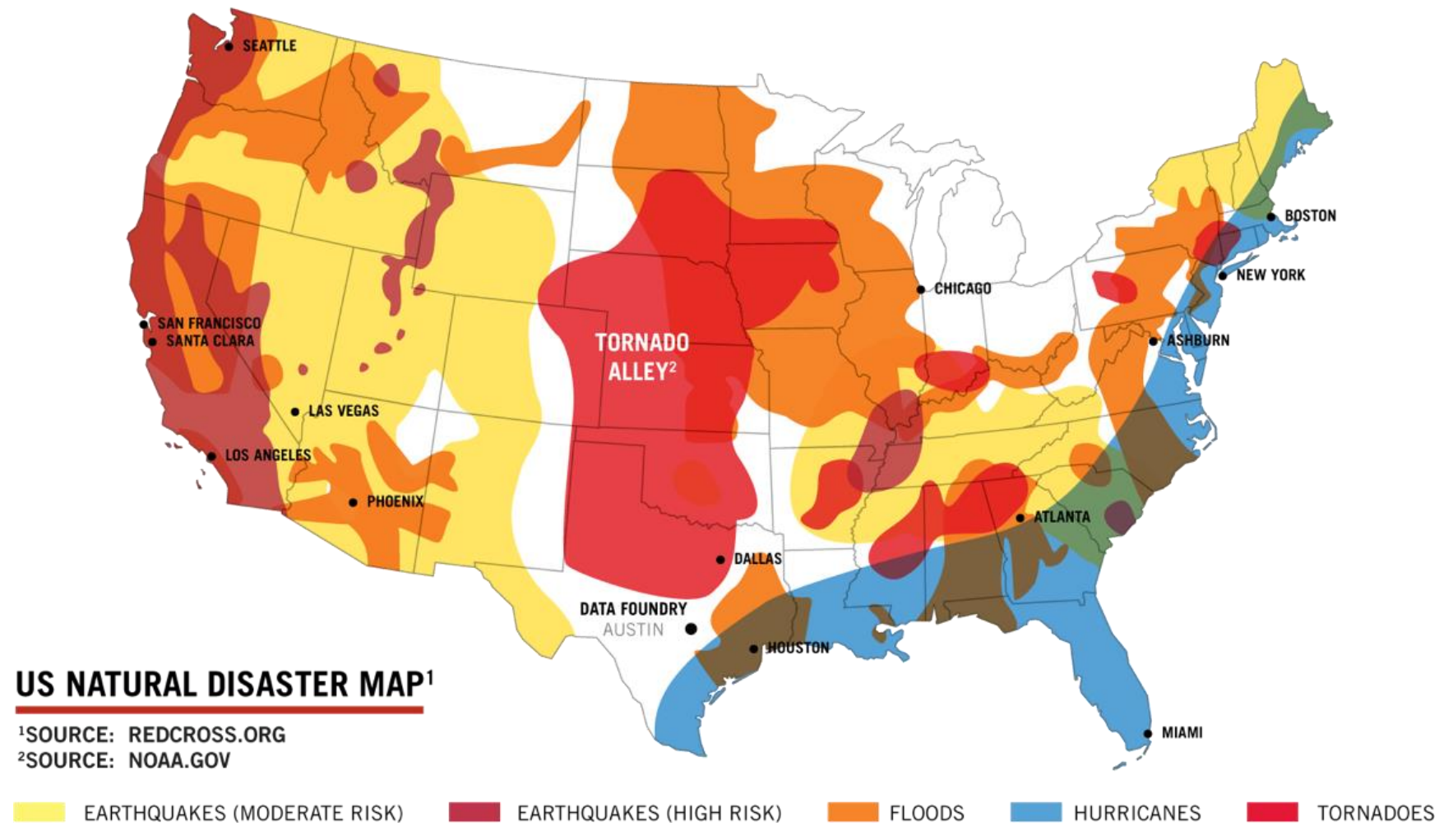

**Figure 4.** U.S. natural disaster risk map [7].

In the experiments, the capacities of the links were assumed to have a uniform distribution, each with a specified range [23,30]. For each disaster scenario, at first link capacities were randomly drawn from their corresponding distributions. Then to replicate the impact of the disaster, the capacities of 50% of the links in the risk areas were further reduced; these links are randomly selected. The reduction in capacity could be as high as 100%, if the objective is to make a link impassable. Since the network employed for the experiments is simplified, there are fewer alternate paths between the O-D pairs. For this reason, an 80% reduction in capacity was assumed to avoid a complete gridlock. Other studies have also used a similar approach for capacity reductions (e.g., [27,30]). The aim of these experiments is to measure at a very high level how the different natural disasters impact freight logistics, for which limited information is available in the literature. Once the general relationship between network performance and disaster occurrence is established, future work can focus on examining specific cases such as comparing the cost of a hurricane in the Gulf Coast (e.g., Hurricane Harvey) versus one in the Southeastern region (e.g., Hurricane Florence) versus one in the Northeastern region (e.g., Hurricane Sandy).

5.2. Results and Discussion

To apply the SAA algorithm, the number of independent disruption-scenario samples ($M$) was set to 100, the sample size ($N$) was set to 1, and the number of large-size samples ($N'$) was set to 1,000. For the GP algorithm, the value of the relative gap (i.e., convergence criterion) was set to 0.0001 [43], which is the relative change in the value of the objective function from one iteration to the next. Note that the terms in the objective function were normalized to yield consistent units. Specifically, the first term was divided by the sum of truck demand and intermodal truck demand and the second term divided by the sum of rail demand and intermodal rail demand.

With the above parameters, the SAA method produced several candidate freight flow patterns, but no more than 100 ($M = 100$). Among these candidate flow patterns, the optimal flow pattern is the one that yields the lowest optimality gap (Equation 16) when each candidate flow pattern is applied to the 1,000 test scenarios ($N' = 1{,}000$).

Table 2 summarizes the cost statistics for the four disaster scenarios. These cost statistics were estimated using samples of different sizes as required by the SAA method. The CPU run times for the four disaster scenarios (high-risk earthquake, high and moderate risk earthquake, hurricane, and flood) were 595.9, 716.2, 669.2, and 417.1 minutes, respectively. As shown, the impact of hurricane is least costly (mean total cost = 48 hours/day) and flood is most costly (mean total cost = 200 hours/day).

**Table 2.** Cost statistics for solutions under different disasters.

| Total Cost (Hour/day) | Earthquake (High Risk) | Earthquake (High and Moderate Risk) | Hurricane | Flood |
|---|---|---|---|---|
| Average | 50.0401 | 76.2006 | 47.9100 | 199.1450 |
| Std. dev. | 0.0579 | 0.1524 | 0.1294 | 0.3699 |
| Minimum | 47.7146 | 70.0737 | 42.7106 | 184.2753 |
| Maximum | 54.4278 | 87.7608 | 57.7205 | 227.2015 |
| gap | 0.2001 | 0.4912 | 0.4162 | 1.1830 |
| $\sigma_{\text{gap}}$ | 0.1939 | 0.5147 | 0.4355 | 1.2484 |

Before comparing the resulting freight flows for road and rail networks for different disaster scenarios, the flows under base condition (i.e., without disaster) are shown in Figure 5. The thickness of the links signifies the volume of assigned freight traffic. The result in left part of Figure 5 indicates that there is high truck flow in Interstates that traverse Arkansas, California, Connecticut, Florida, Georgia, Illinois, Indiana, Michigan, New Jersey, New York, Tennessee, Texas, and Washington. The result in right part of Figure 5 indicates that there is high train flow on rail tracks that traverse Arizona, Georgia, Illinois, Indiana, Iowa, Missouri, Montana, Nebraska, New Jersey, New Mexico, New York, North Dakota, Ohio, Pennsylvania, Texas, Wyoming, and Washington.

The resulting user-equilibrium flow for road and rail networks for different disaster scenarios are shown in Figures 6 through 9. The result in left part of Figure 6 indicates that there is high truck flow on Interstates that traverse Arizona, California, Florida, Georgia, Idaho, Illinois, Indiana, Michigan, New York, Ohio, and Wyoming under the high-risk earthquake scenario. The high truck flow on I-80 in Nevada and Utah is due to freight being diverted from I-5 in California when there is an earthquake. The result in right part of Figure 6 indicates that there is high train flow on rail tracks that traverse Illinois, Indiana, Iowa, Minnesota, Montana, North Dakota, South Dakota, West Virginia, and Wisconsin under the high-risk earthquake scenario. Compared to the base case scenario (without any disaster), there is little difference in the train flow because the rail tracks in these states are not affected by the earthquake in California.

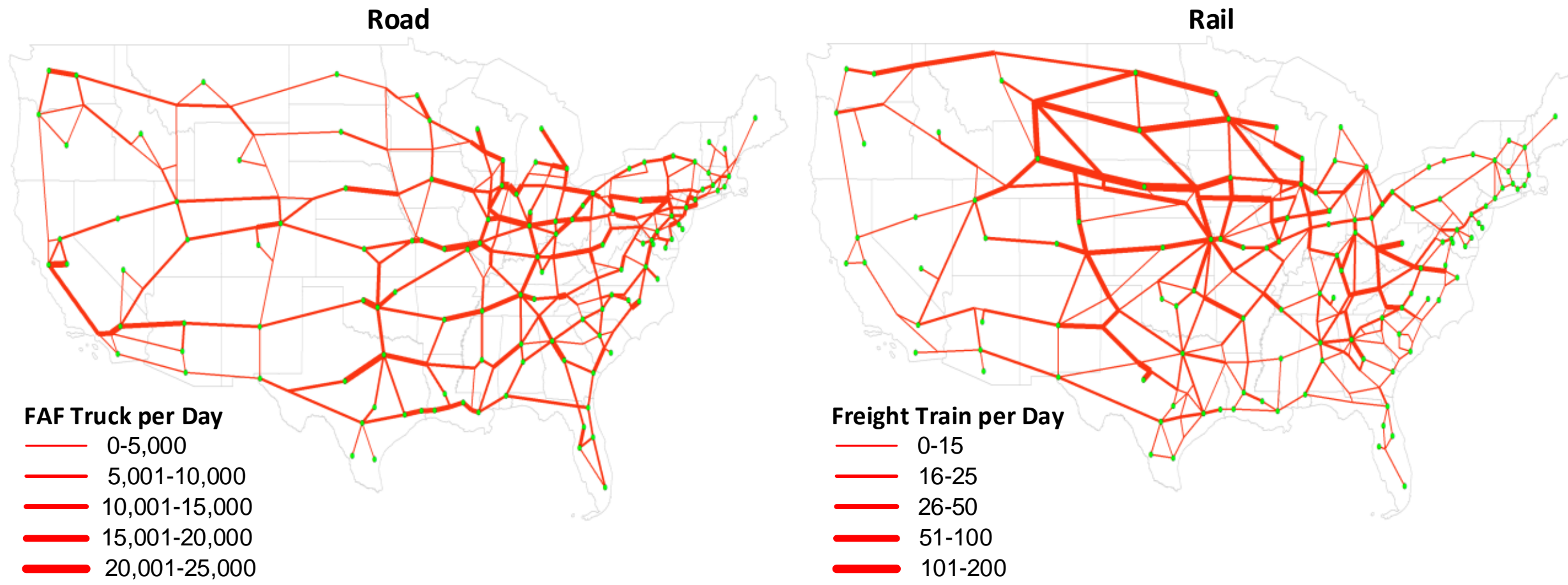

**Figure 5.** Freight traffic assignment (base case).

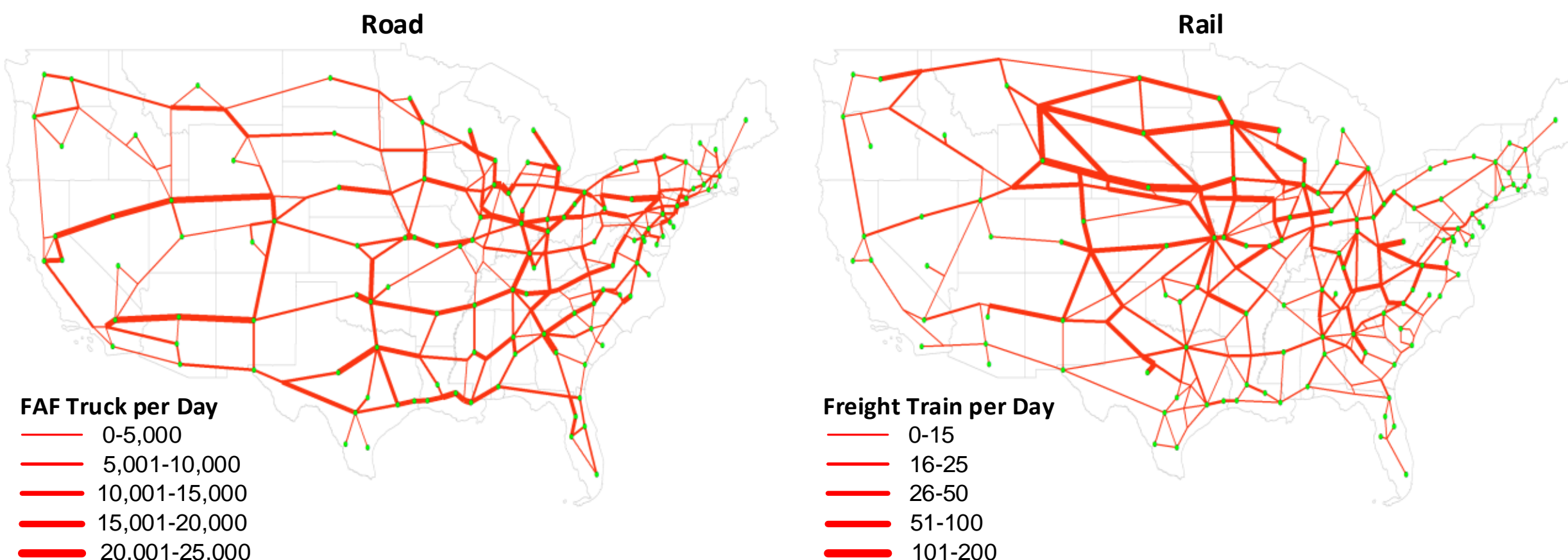

**Figure 6.** Freight traffic assignment under earthquake (high risk).

Figure 7 shows the assigned freight flow under the high and moderate risk earthquake scenario. The result in left part of Figure 7 indicates that there is high truck flow on Interstates that traverse Georgia, Indiana, Kentucky, Pennsylvania, Tennessee, and Texas. Compared to the high-risk earthquake scenario, there is a more even distribution of truck flow in the Western states (such as California, Nevada, Arizona, and Utah). The result in right part of Figure 7 indicates that there is high train flow on rail tracks that traverse Illinois, Iowa, Minnesota, Montana, North Dakota, South Dakota, West Virginia, and Wisconsin. This assigned rail flow is very similar to that of the high-risk earthquake scenario.

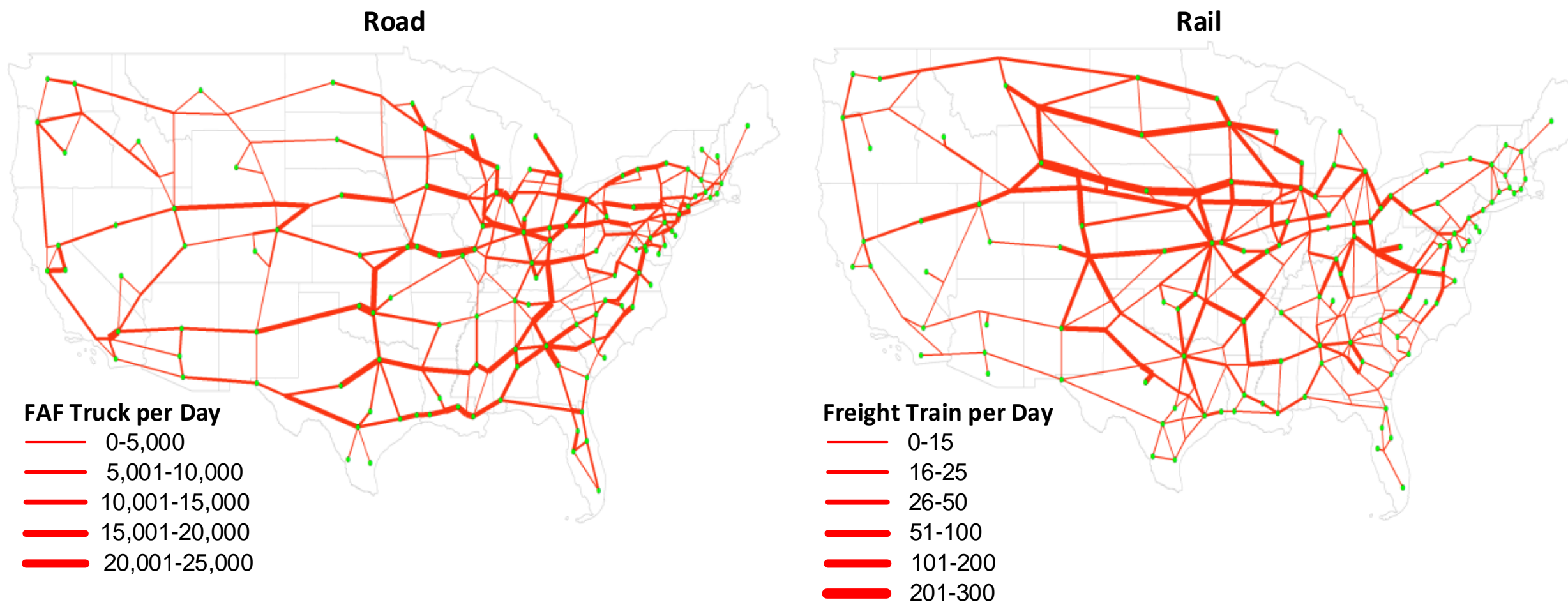

**Figure 7.** Freight traffic assignment under earthquake (high and moderate risk).

Figure 8 shows the assigned freight traffic flow under the hurricane scenario. The result in left part of Figure 8 indicates high truck flow on Interstates that traverse California, Illinois, Indiana, Missouri, Ohio, Tennessee, Texas, and Pennsylvania. Compared to the base case, trucks are diverted from the East and Gulf Coast to the North when there is a hurricane in these regions. The result in right part of Figure 8 indicates high train flow on rail tracks that traverse Iowa, Minnesota, Nebraska, North Dakota, South Dakota, and Wyoming. As is the case with truck flow, there is a higher concentration of rail flow in the Midwest regions.

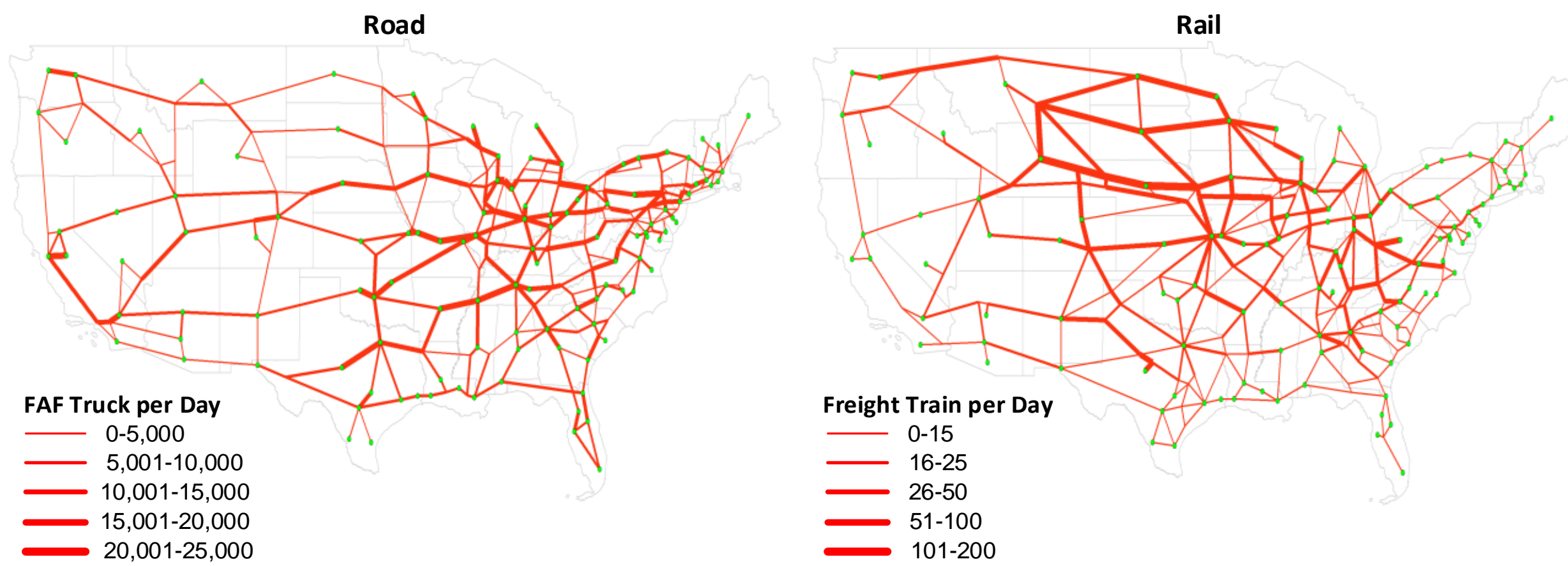

**Figure 8.** Freight traffic assignment under hurricane.

Lastly, Figure 9 shows the assigned freight traffic flow under the flooding scenario. The result in left part of Figure 9 indicates that there is high truck flow in Interstates that traverse Alabama, Arkansas, Indiana, New Mexico, New York, Oklahoma, Pennsylvania, Tennessee, and Texas. Some of the Interstates have very high truck flow; particularly, I-40 in Arkansas and Oklahoma, and I-90 in New York. The reason that trucks are diverting from the Interstates that traverse the Midwestern states is because there is a higher percentage of links in these states that are affected by the flood. The result in right part of Figure 9 indicates that some of the rail tracks have very high train flow (i.e., more than 200 trains per day); particularly, rail tracks in Montana and Wyoming. Furthermore, most of the Mountain states have high rail flow through their states under the flooding scenario. This is also because the trains are avoiding the use of rail tracks in the Midwest regions.

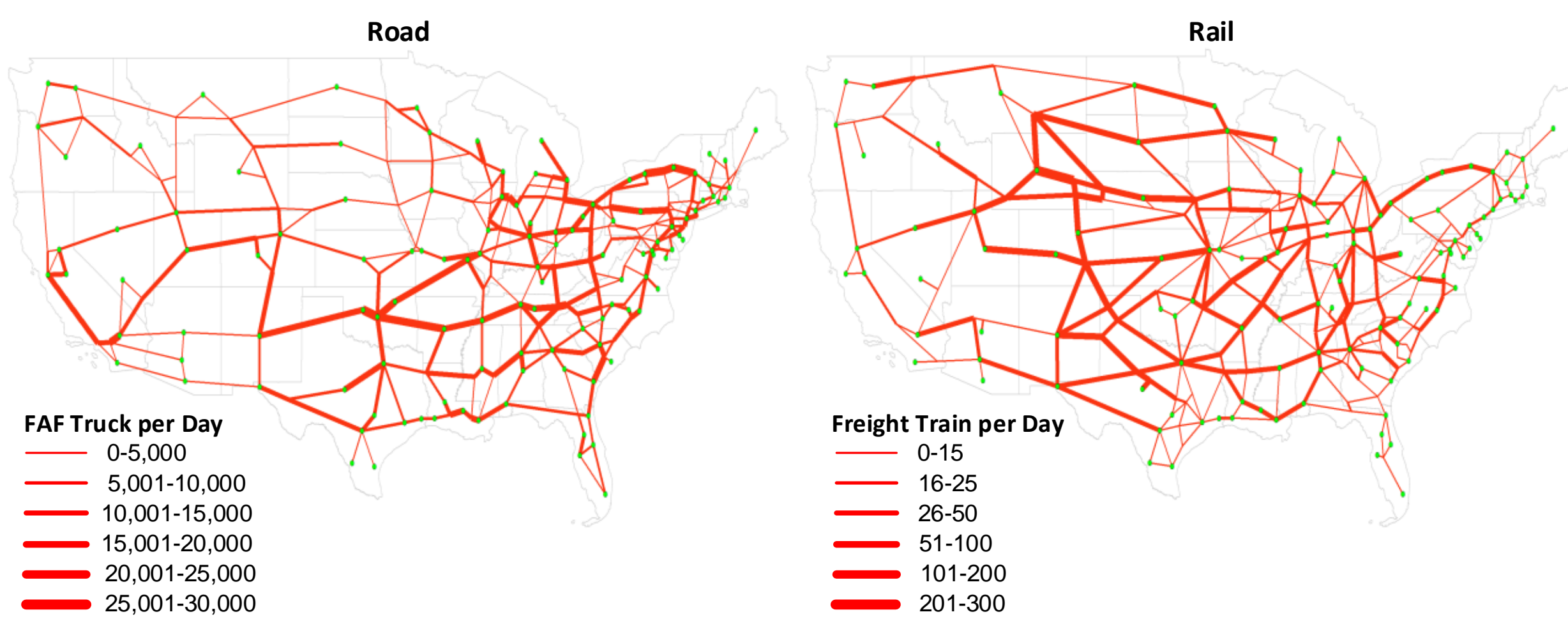

**Figure 9.** Freight traffic assignment under flood.

The proposed model's projected ton-miles under different disaster scenarios are compared quantitatively against those reported in Uddin and Huynh [3], base case in the experiments, for U.S. Census Regions. As evident from Table 3, for both highway and railway modes, the freight ton-miles in case of contiguous U.S. are higher than that of the base case because network uncertainty was not considered. In Midwest region, the highway freight ton-miles is 3% lower for high-risk earthquake, 6% higher for high and moderate risk earthquake, 6% higher for hurricane, and 4% lower for flood compared to that of the base case. In Northeast region, the highway freight ton-miles is 1% higher for high-risk earthquake, 2% higher for high and moderate risk earthquake, 5% higher for hurricane, and 12% higher for flood compared to that of the base case. In South region, the highway freight ton-miles is 8% higher for high-risk earthquake, 7% higher for high and moderate risk earthquake, 2% higher for hurricane, and 25% higher for flood compared to that of the base case. Lastly, in West region, the highway freight ton-miles is 28% higher for high-risk earthquake, 16% higher for high and moderate risk earthquake, 1% higher for hurricane, and 25% higher for flood compared to that of the base case. Overall, for high-risk and high and moderate risk earthquake, West region had the highest increase in freight truck ton-miles. This is expected given the majority of the western states are in earthquake-risk area. For flood, Northeast, South, and West region had a high increase in ton-miles (more than 10%). This is also expected due to the same reason stated above.

**Table 3.** Millions of freight ton-miles under different disasters.

| Mode | Census Region | Uddin and Huynh [3] | Earthquake (High Risk) | Earthquake (High and Moderate Risk) | Hurricane | Flood |
|---|---|---|---|---|---|---|
| Truck[a] | Midwest | 700,458 | 681,225 | 743,699 | 743,790 | 674,983 |
| | Northeast | 224,380 | 225,759 | 227,946 | 235,671 | 251,022 |
| | South | 799,542 | 862,185 | 852,833 | 812,635 | 999,975 |
| | West | 448,319 | 574,544 | 518,351 | 453,333 | 587,796 |
| | Contiguous U.S. | 2,172,699 | 2,343,713 | 2,342,829 | 2,245,429 | 2,513,776 |
| Rail[b] | Midwest | 817,607 | 843,611 | 837,397 | 822,495 | 720,620 |
| | Northeast | 33,328 | 34,271 | 43,041 | 33,940 | 57,010 |
| | South | 417,587 | 422,405 | 435,069 | 429,813 | 595,740 |
| | West | 434,515 | 443,552 | 458,557 | 437,808 | 671,308 |
| | Contiguous U.S. | 1,703,037 | 1,743,839 | 1,774,064 | 1,724,056 | 2,044,678 |

[a]Includes truck, and multiple modes and mail; [b]Includes rail, and multiple modes and mail.

In Midwest region, the rail freight ton-miles is 3% higher for high-risk earthquake, 2% higher for high and moderate risk earthquake, 1% higher for hurricane, and 12% lower for flood compared to that of the base case. In Northeast region, the rail freight ton-miles is 3% higher for high-risk earthquake, 29% higher for high and moderate risk earthquake, 2% higher for hurricane, and 71% higher for flood compared to that of the base case. In South region, the rail freight ton-miles is 1% higher for high-risk earthquake, 4% higher for high and moderate risk earthquake, 3% higher for hurricane, and 43% higher for flood compared to that of the base case. Lastly, in West region, the rail freight ton-miles is 2% higher for high-risk earthquake, 6% higher for high and moderate risk earthquake, 1% higher for hurricane, and 54% higher for flood compared to that of the base case. Overall, the impact of flooding is the highest for Northeast, South, and West regions because there are more states in the flood-risk areas and they are scattered throughout these regions.

## 6. Conclusion

The current freight forecasting methodologies do not consider the risks from weather-induced disruptions which have dramatically increased in recent years; several have occurred recently that severely affected the U.S. freight transportation network. To address these disruptions, this paper developed a stochastic model to assign freight traffic in a large-scale road-rail intermodal network that is subject to network uncertainty (i.e., natural disasters or disruptions). For a specific disaster scenario and given a set of freight demands between origins and destinations and designated modes (road-only, rail-only, and intermodal), the model finds the user-equilibrium freight flow. This paper also provided an algorithmic framework, based on the Sample Average Approximation and Gradient Projection algorithm, to solve the model. Four disaster scenarios were considered in the numerical experiments: high-risk earthquake, high and moderate risk earthquake, hurricane, and flood. The proposed model and algorithmic framework were tested using the U.S. road-rail intermodal network and the Freight Analysis Framework data. The results indicated that when disasters are considered the freight ton-miles are higher than when no disaster is considered, which is expected. The resulting user-equilibrium flows clearly indicate the impact of disasters; that is, truck and rail flow are shifted away from the impacted areas. These results highlight the need to address highways and rail tracks in areas that are normally underutilized but heavily used by trucks and trains when there is a disaster. In terms of cost and freight ton-miles, the impact of flooding is the highest.

This study has a few limitations that could be enhanced in future work. First, the disaster scenarios were generated considering their impact on all intermodal network elements located in the disaster risk area. In reality, a disaster is likely to be more concentrated in one area. Second, the impact of disasters was assumed to be uniform for all randomly selected elements. In reality, the impact will be more severe at the center of the disaster (e.g., eye of a hurricane or epicenter of an earthquake). Lastly, time constraints at the intermodal terminals and/or destinations were not considered. With disaster-related data and freight data becoming more readily available, future work could also compare the freight-related cost of one hurricane versus another; for example, a hurricane that makes landfall in South Florida versus one that makes landfall in the Carolinas.

## Funding

The Authors did not receive any specific funding for this work.

## Data Availability

Freight Analysis Framework data is available from https://www.bts.gov/faf.

## Author Contributions

Conceptualization: M.U. & N.N.H.;
Data curation: M.U.;
Formal analysis: M.U.;
Methodology: M.U.;
Validation: M.U.;
Visualization: M.U.;
Writing – original draft: M.U. & F.A.;
Writing – review & editing: M.U. & N.N.H.

## Conflicts of Interest

The authors have no conflict of interest to declare.

## References

1. U.S. Department of Transportation, Bureau of Transportation Statistics (2022). *Freight Facts and Figures*. https://data.bts.gov/stories/s/Moving-Goods-in-the-United-States/bcyt-rqmu (accessed 1 November 2023).
2. Hwang, T., & Ouyang, Y. (2014). Assignment of freight shipment demand in congested rail networks. *Transportation Research Record*, *2448*(1), 37–44. https://doi.org/10.3141/2448-05
3. Uddin, M.M., & Huynh, N. (2015). Freight traffic assignment methodology for large-scale road-rail intermodal networks. *Transportation Research Record*, *2477*(1), 50–57. https://doi.org/10.3141/2477-06
4. U.S. Department of Transportation, Federal Highway Administration (2023). *Freight Disruptions: Impacts on Freight Movement from Natural and Man-Made Events*. https://ops.fhwa.dot.gov/freight/freight_analysis/fd/index.htm (accessed 1 November 2023).
5. Barton, T. (2018). *Nearly 200 SC roads, including stretch of I-95, closed as flooding worsens*. https://www.thestate.com/news/state/south-carolina/article218717275.html (accessed 1 November 2023).
6. NOAA National Centers for Environmental Information (2023). *U.S. Billion-Dollar Weather and Climate Disasters*. https://www.ncei.noaa.gov/access/billions/summary-stats/US/2022 (accessed 1 November 2023).
7. Disaster Summit (2023). *U.S. Disaster Map*. https://disastersummit.org/us-natural-disaster-map/ (accessed 1 November 2023).
8. Friesz, T.L., Gottfried, J.A., & Morlok, E.K. (1986). A sequential shipper–carrier network model for predicting freight flows. *Transportation Science*, *20*(2), 80–91.
9. Guelat, J., Florian, M., & Crainic, T.G. (1990). A multimode multiproduct network assignment model for strategic planning of freight flows. *Transportation Science*, *24*(1), 25–39.
10. Crainic, T. G., Ferland, J., & Rousseau, J. (1984). A tactical planning model for rail freight transportation. *Transportation Science*, *18*(2), 165–184.
11. Winebrake, J.J., Corbett, J.J., Falzarano, A., Hawker, J.S., Korfmacher, K., Ketha, S., & Zilora, S. (2008). Assessing energy, environmental, and economic tradeoffs in intermodal freight transportation. *Journal of the Air & Waste Management Association*, *58*(8), 1004–1013. https://doi.org/10.3155/1047-3289.58.8.1004
12. Arnold, P., Peeters, D., & Thomas, I. (2004). Modelling a rail/road intermodal transportation system. *Transportation Research Part E: Logistics and Transportation Review*, *40*(3), 255–270. https://doi.org/10.1016/j.tre.2003.08.005
13. Osorio-Mora, A., Núñez-Cerda, F., Gatica, G., & Linfati, R. (2020). Multimodal capacitated hub location problems with multi-commodities: an application in freight transport. *Journal of Advanced Transportation*, *2431763*. https://doi.org/10.1155/2020/2431763
14. Chang, T.-S. (2008). Best routes selection in international intermodal networks. *Computers & Operations Research*, *35*(9), 2877–2891. https://doi.org/10.1016/j.cor.2006.12.025
15. Mahmassani, H., Zhang, K., Dong, J., Lu, C., Arcot, V.C., & Miller-Hooks, E.D. (2007). Dynamic network simulation-assignment platform for multiproduct intermodal freight transportation analysis. *Transportation Research Record*, *2032*(1), 9–16. https://doi.org/10.3141/2032-02
16. Zhang, K., Nair, R., Mahmassani, H., Miller-Hooks, E.D., Arcot, V.C., Kuo, A., Dong, J., & Lu, C. (2008). Application and validation of dynamic freight simulation-assignment model to large-scale intermodal rail network: pan-European case. *Transportation Research Record*, *2066*(1), 9–20. https://doi.org/10.3141/2066-02
17. Loureiro, C.F.G., & Ralston, B. (1996). Investment Selection Model for Multicommodity Multimodal Transportation Networks. *Transportation Research Record*, *1522*(1), 38–46. https://doi.org/0.1177/0361198196152200105

18. Yamada, T., Russ, B.F., Castro, J., & Taniguchi, E. (2009). Designing multimodal freight transport networks: a heuristic approach and applications. *Transportation Science*, *43*(2), 129–143.
19. Fernandez, J.E., De Cea, J., & Giesen, R. (2004). A strategic model of freight operations for rail transportation systems. *Transportation Planning and Technology*, *27*(4), 231–260. https://doi.org/10.1080/0308106042000228743
20. Agrawal, B.B., & Ziliaskopoulos, A. (2006). Shipper–carrier dynamic freight assignment model using a variational inequality approach. *Transportation Research Record*, *1966*(1), 60–70. https://doi.org/10.1177/0361198106196600108
21. Garg, M., & Smith, J.C. (2008). Models and algorithms for the design of survivable multicommodity flow networks with general failure scenarios. *Omega*, *36*(6), 1057–1071. https://doi.org/10.1016/j.omega.2006.05.006
22. Gedik, R., Medal, H., Rainwater, C., Pohl, E.A., & Mason, S.J. (2014). Vulnerability assessment and re-routing of freight trains under disruptions: A coal supply chain network application. *Transportation Research Part E: Logistics and Transportation Review*, *71*, 45–57. https://doi.org/10.1016/j.tre.2014.06.017
23. Uddin, M.M., & Huynh, N. (2016). Routing model for multicommodity freight in an intermodal network under disruptions. *Transportation Research Record*, *2548*(1), 71–80. https://doi.org/10.3141/2548-09
24. Ke, G.Y. (2022). Managing rail-truck intermodal transportation for hazardous materials with random yard disruptions. *Annals of Operations Research*, *309*, 457–483. https://doi.org/10.1007/s10479-020-03699-1
25. Sun, Y., Yu, N., & Huang, B. (2022). Green road–rail intermodal routing problem with improved pickup and delivery services integrating truck departure time planning under uncertainty: An interactive fuzzy programming approach. *Complex & Intelligent Systems*, *8*, 1459–1486. https://doi.org/10.1007/s40747-021-00598-1
26. SteadieSeifi, M., Dellaert, N. P., Nuijten, W., Van Woensel, T., & Raoufi, R. (2014). Multimodal freight transportation planning: A literature review. *European Journal of Operational Research*, *233*(1), 1–15. https://doi.org/10.1016/j.ejor.2013.06.055
27. Chen, L., & Miller-Hooks, E.D. (2012). Resilience: an indicator of recovery capability in intermodal freight transport. *Transportation Science*, *46*(1), 109–123. https://doi.org/10.1287/trsc.1110.0376
28. Misra, S., Padgett, J.E. (2019). Seismic resilience of a rail-truck intermodal freight network. *13th International Conference on Applications of Statistics and Probability in Civil Engineering*.
29. Huang, M., Hu, X., & Zhang, L. (2011). A decision method for disruption management problems in intermodal freight transport. In J. Watada, G. Phillips-Wren, L.C. Jain, R.J. Howlett (Eds.), *Intelligent Decision Technologies*, Springer, Berlin Heidelberg.
30. Miller-Hooks, E.D., Zhang, X., & Faturechi, R. (2012). Measuring and maximizing resilience of freight transportation networks. *Computers & Operations Research*, *39*(7), 1633–1643. https://doi.org/10.1016/j.cor.2011.09.017
31. Uddin, M., & Huynh, N. (2019). Reliable routing of road–rail intermodal freight under uncertainty. *Networks and Spatial Economics*, *19*, 929–952. https://doi.org/10.1007/s11067-018-9438-6
32. Kramarz, M., Przybylska, E., & Wolny, M. (2022). Reliability of the intermodal transport network under disrupted conditions in the rail freight transport. *Research in Transportation Business & Management*, *44*, 100686. https://doi.org/10.1016/j.rtbm.2021.100686
33. Ke, G. Y. (2022). Managing rail-truck intermodal transportation for hazardous materials with random yard disruptions. Annals of Operations Research, 309(2), 457-483. https://doi.org/10.1007/s10479-020-03699-1
34. Misra, S., & Padgett, J. E. (2022). Estimating extreme event resilience of rail–truck intermodal freight networks: Methods, models, and case study application. *Journal of Infrastructure Systems*, 28(2), 04022013. https://doi.org/10.1061/(ASCE)IS.1943-555X.0000679
35. Tafur, A., Padgett, J. E., & Gonzalez, A. Seismic Resilience of Intermodal Freight Transportation Networks with Ties to Economic Impact Modeling. In Lifelines 2022 (pp. 372-383).
36. Sheffi, Y. (1985). *Urban Transportation Networks: Equilibrium Analysis with Mathematical Programming Methods*, Prentice Hall, Inc., Upper Saddle River, N.J.
37. Santoso, T., Ahmed, S., Goetschalckx, M., & Shapiro, A. (2005). A stochastic programming approach for supply chain network design under uncertainty. *European Journal of Operational Research*, *167*(1), 96–115. https://doi.org/10.1016/j.ejor.2004.01.046
38. Shapiro, A. (1996). Simulation-based optimization—convergence analysis and statistical inference. *Communications in Statistics. Stochastic Models*, *12*(3), 425–454. https://doi.org/10.1080/15326349608807393
39. Jayakrishnan, R., Tsai, W.K., Prashker, J.N., & Rajadhyaksha, S. (1994). Faster path-based algorithm for traffic assignment. *Transportation Research Record*, *1443*(1), 75–83.
40. Center for Transportation Analysis, Oak Ridge National Laboratory. (2018). *Intermodal transportation network*. https://www-cta.ornl.gov/transnet/Intermodal_Network.html (accessed 27 May 2018).
41. U.S. Department of Transportation, Federal Highway Administration. (2007). Freight Analysis Framework, FAF3 [datasets].
42. Uddin, M., & Huynh, N. (2020). Model for collaboration among carriers to reduce empty container truck trips. *Information*, *11*(8), 377. https://doi.org/10.3390/info11080377
43. Boyce, D., Ralevic-Dekic, B., & Bar-Gera, H. (2004). Convergence of traffic assignments: how much is enough? *Journal of Transportation Engineering*, *130*(1), 49–55. https://doi.org/10.1061/(ASCE)0733-947X(2004)130:1(49)